\documentclass[10pt]{article}

\def\Frac#1#2{\displaystyle{\frac{#1}{#2}}}

\newtheorem{theorem}{Theorem}

\newtheorem{remark}{Remark}
\usepackage{amssymb}

\begin{document}

\title{A new type of sharp bounds for ratios of modified Bessel functions}

\author{Diego Ruiz-Antol\'{\i}n\\
Departamento de Matem\'aticas, Estad\'{\i}stica y 
        Computaci\'on,\\
        Univ. de Cantabria, 39005 Santander, Spain.\\ 
\and
Javier Segura\\
        Departamento de Matem\'aticas, Estad\'{\i}stica y 
        Computaci\'on,\\
        Univ. de Cantabria, 39005 Santander, Spain.\\   
{ \small
  e-mail: {\tt
    diego.ruizantolin@unican.es,
    javier.segura@unican.es}}
    }

\date{\today}

\maketitle
\begin{abstract}
The bounds for the ratios of first and second kind modified Bessel functions of consecutive 
orders are important quantities appearing in a large number of scientific applications.
We obtain new bounds which are accurate in a large region
of parameters and which are shaper than previous bounds. The new bounds are obtained by a qualitative analysis of the
Riccati equation satisfied by these ratios. A procedure is considered in
which the bounds obtained from the analysis of the Riccati equation are used to
define a new function satisfying a new Riccati equation which yields sharper 
bounds. Similar ideas can be applied to other functions.
\end{abstract}

\vskip 0.8cm \noindent
{\small
\par\noindent
Keywords \& Phrases:
Modified Bessel functions; bounds.
}

\section*{Introduction}

The ratios of modified Bessel functions $I_{\nu}(x)/I_{\nu-1}(x)$ and
$K_{\nu}(x)/K_{\nu-1}(x)$ are important quantities appearing in a large number of scientific applications. 
Bounds for these ratios have been recently used 
in connection with Schwarz methods for reaction-diffusion processes \cite{Gigante:2013:OSM}, statistics \cite{Igarashi:2014:ROI,Pal:2014:GEF}
and in the study of oscillatory solutions of second order ODEs \cite{Heitman:2015:OTE}. 
See \cite{Segura:2011:BFR} and references
cited therein for additional examples of application of these bounds; see also \cite{Gronwall:1932:AIF, Amos:1974:COM, Nasell:1978:RBF, Laforgia:2010:SIF, Hornik:2013:ATB, Baricz:2015:BFT} for additional papers exploring different methods 
for bounding these ratios.

In \cite{Segura:2011:BFR}, several techniques were considered for bounding the ratios of modified Bessel functions. 
One of the most useful ideas was the use of the Riccati equations satisfied by these ratios and, in particular, the 
fact that one of the critical curves given by $h'(x)=0$ (with $h(x)$ the ratio under consideration) is a bound for the ratio.
In this paper, this use of the Riccati equation is extended and improved in two different ways. 

Firstly, we will
show how this technique alone is able to provide both sharp lower and upper bounds (only one-sided inequalities were obtained
with this technique in \cite{Segura:2011:BFR}). These bounds are described by a uniparametric family of functions, and the best
possible upper and lower bounds in this family are obtained. In the second place, it will be shown how the bounds obtained in the first stage
can be used, after a change of function, to obtain a new Riccati equation which provides a new and improved uniparametric family. 
The bounds obtained in the first stage are also members of the second
uniparametric family.

For the first kind Bessel function, the best possible bounds of the 
second family are sharper than the bounds obtained in the first stage, and particularly in the
limits $x\rightarrow 0,+\infty$ and $\nu\rightarrow +\infty$. For the second kind Bessel function the situation
is similar, with the difference that only the upper bound can be improved in the second stage. 

\section{Bounds from Riccati equations}
\label{bori}

Bounds for the ratios of modified Bessel functions can be obtained by considering the following result (which we prove
for completeness):

\begin{theorem}
\label{esode}
Let $h(x)$ be a solution of the Riccati 
equation $h'(x)=A(x)+B(x)h(x)+C(x)h(x)^2$ defined in $(0,+\infty)$, where the coefficients are continuous and satisfying $A(x)C(x)<0$. 
Let $\phi (x)$ be the positive root of $A(x)+B(x)\phi (x)+C(x)\phi(x)^2=0$; if $\phi (x)$
is strictly monotonic the following holds:
\begin{enumerate}
\item{If $C(x)<0$}, $h(0^+)>0$ and $\phi' (0^+) h ' (0^+)>0$, then $h(x)<\phi(x)$ if $\phi '(x)>0$ and
  $h(x)>\phi(x)$ if $\phi '(x)<0$
\item{If $C(x)>0$}, $h(+\infty)>0$ and $\phi' (+\infty) h ' (+\infty)>0$, then $h(x)<\phi(x)$ if $\phi '(x)<0$ and
  $h(x)>\phi(x)$ if $\phi '(x)>0$
\end{enumerate} 
\end{theorem}

\noindent {\it Proof.} We prove the case $C(x)<0$ and $\phi '(x)>0$, the remaining three cases are proved
in the same way. 

Because $C(x)<0$ the fact that $\phi' (0^+) h ' (0^+)>0$ (and then $h'(0^+)>0$) implies that
$0<h(0^+)<\phi(0^+)$, but then, because $\phi (x)$ is increasing,
this must be true for all $x>0$, and therefore $h(x)<\phi(x)$ for positive $x$. Indeed, a value $x_0$ such that $h (x_0) = \phi (x_0)$ 
(and therefore $h'(x_0)=0$) 
can not be reached, because $h (x)$ would approach the graph of the increasing function $\phi (x)$ from below and for this it is 
necessary that $h' (x_0) > \phi ' (x_0) >0$, in contradiction with the fact that $h '(x_0)=0$.
\hfill $\square$.

It is not essential that $A(x)C(x)<0$ for this type of result to hold; we refer to \cite{Segura:2012:OBF} for a more general result. 
However, this simple theorem is all that is required to find a first iteration of upper and lower bounds. In particular, as we will see,
considering the function $h(x)=x^{-\alpha}I_{\nu}(x)/I_{\nu-1}(x)$ for some values of $\alpha$ we will re-obtain the best possible upper
and lower bounds of the form 
$x/(\lambda+\sqrt{\lambda^2+x^2})$ 
for $I_{\nu}(x)/I_{\nu-1}(x)$, and similarly for $K_{\nu}(x)/K_{\nu-1}(x)$.
These are sharp bounds as $x$ and $\nu$ become large, and also in some cases when $x$ goes to zero for the $I$ function and as $\nu$ approaches $1/2$
for the $K$ function.

This type of bounds is the first step in the generation of new and improved bounds. The idea will be to use these first bounds for 
redefining the function and to analyze the new Riccati equation. We explain this procedure next.

\subsection{Iterating the bounds from Riccati equations}
\label{iterating}

We start from the Riccati equation
$$
h_0 '(x)=A_0 (x)+B_0(x)h_0 (x)+C_0(x)h_0 (x)^2.
$$
Next we consider the function $h_1 (x)=h_0 (x)/\phi_0 (x)$, where $\phi_0 (x)\equiv \beta_0 (x)$ 
is a function of convenience which we 
will choose as one of the roots of $A_0 (x)+B_0(x)\phi_0 (x)+C_0(x)\phi_0 (x)^2=0$ (which will be a bound for $h_0(x)$). 
The next step
is
\begin{equation}
\begin{array}{l}
h_1 '(x)=A_1(x)+B_1(x)h_1 (x)+C_1(x) h_1(x)^2\\
A_1 (x)=\Frac{A_0(x)}{\phi_0 (x)},\,B_1 (x) =B_0(x)-\Frac{\phi_0^{\prime}(x)}{\phi_0(x)},\,C_1 (x)=\phi_0(x) C_0(x)
\end{array}
\end{equation}
If one of the roots $\phi_1 (x)$, solution of characteristic equation $A_1 (x)+B_1 (x)\phi_1 (x)+C_1 (x)\phi_0 (x)^2=0$,
turns out to be a bound for $h_1(x)$, then $\beta_1 (x)=\phi_1 (x)\beta_0 (x)=\phi_1(x)\phi_0(x)$ will be a bound for
$h_1 (x)$.

In general, we could consider the iteration of this process and, after $n$ steps,
\begin{equation}
\label{rico}
\begin{array}{l}
h_n ' (x)=A_n(x)+B_n(x)h_n (x)+C_n(x) h_n (x)^2\\
A_n (x)=\Frac{A_0(x)}{\beta_{n-1} (x)},\,B_n (x) =B_0(x)-\sum_{i=0}^{n-1}\Frac{\phi_i^{\prime}(x)}{\phi_i(x)},\,C_n (x)=\beta_{n-1} (x) C_0(x)
\end{array}
\end{equation}
with
\begin{equation}
\beta_{n-1} (x)=\prod_{i=0}^{n-1}\phi_i (x).
\end{equation}
Then the characteristic roots of (\ref{rico}) can be written as
\begin{equation}
\label{care}
\phi_n (x)=\Frac{2 A_n (x)}{-B_n (x)\pm \sqrt{B_n (x)^2-4A_n (x) C_n (x)}}
\end{equation}
and we have
\begin{equation}
\label{betan}
\beta_n (x)=\beta_{n-1} (x) \phi_n (x)=\Frac{1}{\eta_n (x) \pm \sqrt{\eta_n (x) ^2+\gamma (x)}}
\end{equation}
where
\begin{equation}
\gamma (x)=-C_0(x)/A_0(x)
\end{equation}
and
\begin{equation}
\eta_n (x)=-\Frac{B_n (x)}{2 A_0 (x)}=
\eta_0 (x) \left(1-\Frac{1}{B_0 (x)}\Frac{d}{dx}\log\beta_{n-1} (x)\right).
\end{equation}
This can be written as
\begin{equation}
\label{etar}
\eta_n (x)=\eta_0 (x)\mp\Frac{\eta_{n-1}'(x)+\Frac{1}{2}\beta_{n-1}(x)\gamma'(x)}{2A_0 (x)\sqrt{\eta_{n-1}(x)^2+\gamma (x)}}
\end{equation}
where the signs in (\ref{etar}) are those in correspondence with the previous bound
\begin{equation}
\beta_{n-1} (x)=\Frac{1}{\eta_{n-1}(x)\pm \sqrt{\eta_{n-1}(x)^2+\gamma(x)}}.
\end{equation}

For determining whether the successive values $\beta_n (x)$ are bounds we can invoke Theorem \ref{esode}.
For proving that $\beta_n (x)$ is a bound for $h_n(x)$ we must check the monotonicity of
the characteristic root $\phi_n(x)=\beta_{n}(x)/\beta_{n-1}(x)$, which should correspond to the
 monotonicity property of $h_n(x)=h_0(x)/\beta_{n-1}(x)$ (for $h_n(x)$ we only need to check the 
monotonicity either for $x=0^+$
or $x=+\infty$).

\section{Bounds for the ratio $I_{\nu}(x)/I_{\nu-1}(x)$}

We start with  the Riccati equation satisfied by the
ratio $h (x)=I_{\nu}(x)/I_{\nu -1}(x)$
\begin{equation}
\label{RI}
h^{\prime}(x)=1-\Frac{2(\nu-1/2)}{x}h(x)-h(x)^2,
\end{equation}
which can be easily proved using \cite[10.29.2]{Olver:2010:BF}. 

Because $h(x)$ is positive, the relevant characteristic root of the Riccati equation is
the positive root and it was proved in \cite{Segura:2011:BFR} that it is an upper bound
for the ratio. We can consider a more general situation by starting from
\begin{equation}
\label{h01}
h_0(x)=x^{-\alpha}h (x), 
\end{equation}
which leads to
\begin{equation}
\label{ricaa}
h_0 '(x)=x^{-\alpha}-\Frac{2\lambda}{x}h_0 (x)-x^{\alpha} h_0 (x)^2,\,\lambda=\nu+\frac12 (\alpha -1)
\end{equation}
For this new  equation, and using the notation of section \ref{iterating},
\begin{equation}
\label{otrosc}
\gamma (x)=x^{2\alpha},\,\eta_0 (x) =\lambda x^{\alpha-1}
\end{equation}
and the potential bound for 
$h_0(x)$ is \footnote{Observe that we can not use
the equations in section \ref{iterating} 
for going from the equation for $h (x)$ to that for $h_0 (x)$ because $x^{\alpha}$
is not a characteristic root for the first equation, but we can use these formulas once we obtain the equation for $h_0 (x)$}
\begin{equation}
\label{bb}
\beta_0 (x)=\phi_0 (x)=\Frac{1}{\eta_0(x)+\sqrt{\eta_0 (x)^2+\gamma (x)}}=\Frac{x^{1-\alpha}}{\lambda+\sqrt{\lambda^2 +x^2}}
\end{equation}
and we have a family of potential bounds for $h(x)$ depending on the  parameter $\alpha$:
\begin{equation}
\label{ba}
b_{\alpha} (\nu,x)=\Frac{x}{\lambda+\sqrt{\lambda^2 +x^2}},\,\lambda=\nu+\frac12 (\alpha -1).
\end{equation}
It is known that $b_0 (\nu,x)$ is an upper bound for $h(x)$ when $\nu\ge 1/2$, and this was proved 
\cite{Segura:2011:BFR} using the Riccati equation. On the other hand, $b_{1} (\nu, x)$ is known to be a lower bound as was
proved using different techniques
(see for instance \cite{Laforgia:2010:SIF,Segura:2011:BFR}). We will prove very easily both bounds using Theorem \ref{esode} and
show that they are the best possible bounds of the form (\ref{ba}).

\begin{theorem}
\label{boundeb}
$$
\Frac{I_{\nu}(x)}{I_{\nu-1}(x)}<b_0 (\nu,x)=\Frac{x}{\nu-1/2+\sqrt{(\nu-1/2)^2+x^2}},\,\nu\ge 1/2,
$$
$$
\Frac{I_{\nu}(x)}{I_{\nu-1}(x)}>b_1 (\nu,x)=\Frac{x}{\nu+\sqrt{\nu^2+x^2}},\,\nu\ge 0.
$$
These bounds are the best possible of the form $b_{\alpha} (\nu,x)$, $\alpha\in {\mathbb R}$, 
in their range of validity.
\end{theorem}
\noindent {\it Proof.} 
The first bound corresponds to the case $\alpha=0$ and the second to $\alpha=1$.

We start by noticing that $b_{\alpha}(\nu,x)$ are not bounds if $\alpha\in (0,1)$. 
This is easy to check by considering the Maclaurin series \cite[10.25.2]{Olver:2010:BF} and the asymptotic series 
\cite[10.40.1]{Olver:2010:BF}, leading to
\begin{equation}
\label{h0i}
\begin{array}{l}
h_0 (x)=x^{-\alpha}\left(\Frac{x}{2\nu}-\frac18 \Frac{x^3}{\nu^2 (\nu+1)}+{\cal O}(x^5)\right),\\
\\
h_0 (x)=x^{-\alpha}\left(1-\Frac{\nu-1/2}{x}+{\cal O}(x^{-2})\right)
\end{array}
\end{equation}
From this we see that $h_0'(0^+)h_0'(+\infty)<0$ if $\alpha \in (0,1)$, and this means that the graph of $h_0(x)$ must cross
the graph of the characteristic root $\phi_0(x)$, which therefore can not be a bound for $h_0 (x)$.

Next, we prove the bounds for $\alpha=0$ (upper bound) and $\alpha=1$ (lower bound)

\begin{enumerate}
\item{$\alpha=0$:} the case $\nu=1/2$ is obvious using \cite[10.27.2]{Olver:2010:BF}. For $\nu>1/2$ 
we have that $\phi_0 (x)=x/(\lambda+\sqrt{\lambda^2+x^2})$, $\lambda=\nu -1/2$, is strictly increasing\footnote{Observe that 
this is true only if $\lambda>0$, that is, if $\nu>1/2$}; in
addition, (\ref{h0i}) shows that $h(0^+)>0$ and $h_0 '(0^+)>0$, and therefore $\phi_0 ' (0^+)h_0 '(0^+)>0$. Applying Thm. \ref{esode}, 
we have $h_0(x)=I_{\nu}(x)/I_{\nu-1}(x)<\phi_0 (x)$.
\item{$\alpha=1$:} the case $\nu=0$ is obvious because $I_0(x)/I_{-1}(x)=I_0(x)/I_1(x)$ and the previous inequality
we have proven shows that this ratio is greater than $1$. We observe that $\phi_0 (x)=1/(\nu+\sqrt{\nu^2 +x^2})$
is strictly decreasing as a function of $x$ and that for $\nu>0$ 
(\ref{h0i}) shows that $h(0^+)>0$, $h_0 '(0^+)<0$, and therefore $\phi_0 ' (0^+)h_0 '(0^+)>0$. Applying Thm. \ref{esode}, 
we have $h_0(x)= x^{-1} I_{\nu}(x)/I_{\nu-1}(x)>\phi_0 (x)$.
\end{enumerate}

Finally, we notice that $b_{\alpha} (\nu, x)$ is decreasing as
a function of $\alpha$, and then, because $b_1 (\nu,x)$ is a lower bound, $b_{\alpha} (\nu,x)$ are lower bounds for 
any $\alpha\ge 1$; similarly, because $b_0 (\nu,x)$ is an upper bound, $b_{\alpha} (\nu,x)$ are upper bounds for any 
$\alpha\le 0$. The sharpest bounds correspond to $\alpha=0,1$. 

\hfill $\square$.

We observe that Theorem \ref{esode} can not be used to prove that any of the values $\beta_0 (x)$ is a bound for values
$\alpha \in (-1,0)\cup(0,1)$, because $\phi_0 (x)=\beta_0 (x)$ is not monotonic in this case (the derivative of $\beta_0 (\lambda,x)$ is zero at $x=\pm \sqrt{1-\alpha ^2}\lambda/\alpha$). However, from the monotonicity of $b_\alpha (\nu,x)$ we concluded
that for $\alpha\in (-1,0)$ we have upper bounds for $I_{\nu}(x)/I_{\nu-1}(x)$ 
at least for $\nu\ge 1/2$; in fact, these bounds, although less sharp than 
$b_0(\nu,x)$, extend the range of validity with respect to $\nu$, and for $\alpha=-1$ we obtain the following result:

\begin{theorem}
\label{masdeb}
$$
\Frac{I_{\nu}(x)}{I_{\nu-1}(x)}<b_{-1} (\nu,x)=\Frac{x}{\nu-1+\sqrt{(\nu-1)^2+x^2}},\,\nu\ge 0,
$$
\end{theorem}

\noindent {\it Proof.} For $\nu=0$ the bound holds because 
it is equivalent to $I_1 (x)/I_0 (x)>x/(1+\sqrt{1+x^2})$, which is true due to Theorem  \ref{boundeb} (taking $\nu=1$ in the 
second bound). 
For $\nu>0$ we have that $\phi_0(x)=x b_{-1}(\nu,x)$ is increasing as a function of $x$ and $h_0 '(0^+)>0$; 
Theorem \ref{esode} implies that $h_0(x)=xI_{\nu}(x)/I_{\nu}(x)<\phi_0(x)=x^2/(\nu-1+\sqrt{(\nu-1)^2+x^2})$ \hfill $\square$

\vspace*{0.3cm}
The three-term recurrence relation can be used to obtain sequences of converging bounds, just by writing this as

\begin{equation}
\label{CF}
\Frac{I_{\nu}(x)}{I_{\nu -1}(x)}=\Frac{1}{\Frac{2\nu}{x}+\Frac{I_{\nu+1}(x)}{I_{\nu }(x)}}.
\end{equation}

Substituting $I_{\nu+1}(x)/I_{\nu }(x)$ by one of the upper (lower) bounds we obtain a lower (upper) bound for
$I_{\nu}(x)/I_{\nu-1}(x)$. This process can be continued for obtaining a sequence of convergent bounds,
as described in \cite{Segura:2011:BFR}; the sequence is convergent as a consequence of the fact that the continued
fraction resulting from the iteration of (\ref{CF}) is convergent. 
In \cite{Amos:1974:COM}, these sequences of convergent bounds are also considered
(but formulated in a different way).

For instance, iterating only once and starting with Theorem \ref{boundeb} we have:
\begin{theorem}
\label{boundeb2}
$$
\Frac{x}{\nu-1/2+\sqrt{(\nu+1/2)^2+x^2}}<\Frac{I_{\nu}(x)}{I_{\nu-1}(x)}<\Frac{x}{\nu-1+\sqrt{(\nu+1)^2+x^2}},\,\nu\ge 0.
$$
\end{theorem}
The lower bound in this theorem is an improvement over the bound in Theorem \ref{boundeb} for all $x>0$, $\nu\ge 0$. The
upper bound is an improvement over Theorem \ref{boundeb} only if  $x^2<4\nu(2\nu+1)$. Of course, as more iterations of
the recurrence are considered the bound should improve; however, the convergence is slower for larger $x$.

\subsection{Bounds from the iteration of the Riccati equation}

Next we consider the iteration of the Riccati equation, as described in section \ref{bori}. As we will see, the first
iteration gives bounds that not only are superior to the first iteration of the recurrence, but also superior to any number 
of iterations of the recurrence for large enough $x$. We will obtain a new uniparamentric family of bounds which includes the 
bounds in Theorems \ref{boundeb} and \ref{masdeb}; the best possible bounds in this new family will improve the bounds
of those two previous theorems.

We start from (\ref{h01}) and (\ref{otrosc}). The first iteration (following the notation of section \ref{iterating}) gives
\begin{equation}
\label{eta1I}
\eta_1 (x)=x^{\alpha -1}\left[\lambda\left(1+\Frac{1}{2\sqrt{\lambda^2 +x^2}}\right)-\Frac{\alpha}{2}\right]
\end{equation}
and the next potential bounds for $h_0 (x)$ are $\beta_1 (x)=1/(\eta_1 (x)+\sqrt{\eta_1 (x)^2+x^{2\alpha}})$; this means that
the uniparametric family of potential bounds for $I_{\nu}(x)/I_{\nu-1}(x)=x^{\alpha}h_0(x)$ has the form
\begin{equation}
\begin{array}{l}
B_{\alpha} (\nu,x)=\Frac{x}{\delta_{\alpha} (\nu,x)+
\sqrt{\delta_{\alpha} (\nu ,x)^2+x^2}},\\
\delta_{\alpha} (\nu,x)=(\nu-1/2)+\Frac{\lambda}{2\sqrt{\lambda^2+x^2}},\,\lambda=\nu+(\alpha-1)/2.
\end{array}
\end{equation}
The next theorem shows that these are indeed bounds for $\alpha\notin (0,2)$ and that the best possible bounds 
correspond to $\alpha=0,2$.

\begin{theorem} 
\label{cotasric}
$$
\Frac{I_{\nu}(x)}{I_{\nu-1}(x)}<B_{0}(\nu,x),\,\nu\ge 1/2
$$
$$
\Frac{I_{\nu}(x)}{I_{\nu-1}(x)}>B_{2}(\nu,x),\,\nu\ge 0
$$
in their range of validity. These bounds are the best possible of the form $B_{\alpha}(\nu ,x)$ for fixed
 $\alpha\in {\mathbb R}$, 
\end{theorem}

\noindent {\it Proof.} The first bound corresponds to $\alpha=0$ and the second to $\alpha=2$.

We observe that for $\alpha\in (0,2)$, $B_{\alpha}(\nu,x)$ can not be bounds for all $x>0$ when $\nu$ is sufficiently large.
Indeed, considering series expansions for $h_1(x)=h_0(x)/\beta_0(x)$ (see (\ref{h01}) and (\ref{bb})) we have
\begin{equation}
\label{serh1}
\begin{array}{l}
h_1 (x)=\Frac{\lambda}{\nu}-\Frac{1}{8}\Frac{4\nu (\alpha-2)+(\alpha-1)^2}{\nu^2 (\nu+1)(2\nu+\alpha-1)}x^2+{\cal O}(x^4),\\
\\
h_1(x)=1+\Frac{\alpha}{2}x^{-1}+\Frac{1}{8}\left(\alpha^2+2-4\nu\right)x^{-2}+{\cal O}(x^{-3}),
\end{array}
\end{equation}
where in the first expansion we are assuming that $\lambda\ge 0$.
From this we observe that $h_{1}'(+\infty)<0$ if $\alpha>0$ while $h_{1}'(0^+)>0$ if $\nu>(\alpha-1)^2/(8-4\alpha)>0$; 
in this case $\phi_1 (x)$ can not be a bound for $h_1(x)$ because the graph of $h_1 (x)$ necessarily crosses the graph of
$\phi_1 (x)$ for some positive $x$. 

If we can prove that for $\alpha=0$  we have an upper bound and for $\alpha=2$ a lower bound the theorem will be proved.
The fact that these will be the best possible bounds of the form $B_{\alpha} (\nu,x)$ will be a consequence
of the fact that $B_{\alpha} (\nu,x)$ is decreasing as a function of $\alpha$. Let us then prove that $\alpha=0$ and $\alpha=2$
correspond to an upper and a lower bound respectively.

\vspace*{0.2cm}
\noindent
{\bf 1. {\boldmath $\alpha=0$:}}

The case $\nu=1/2$ holds trivially because $I_{1/2}(x)/I_{-1/2}(x)<1$. Let us now assume that $\nu> 1/2$
and therefore $\lambda>0$. In this case we have $h_1 (0^+)>0$ and $h_1 '(0^+ )>0$, and then, on account
of Theorem \ref{esode}, if $\phi_1^{\prime}(x)>0$ for $x>0$, we will have $h_1 (x)<\phi_1 (x)$ (and then 
$I_{\nu}(x)/I_{\nu-1}(x)=h_0(x)<\phi_0(x)\phi_1(x)=\beta_1 (x)=B_0(\nu ,x)$). All we have to prove is that $\phi_1 '(x)>0$ for $\nu>1/2$.

We write $\phi_1 (x)=\beta_1 (x)/\beta_0 (x)$
and then
$$
\phi_1'(x)=\left[\Frac{\sqrt{\eta_0 (x)^2+1}+\eta_0 (x)}{\sqrt{\eta_1 (x)^2+1}+\eta_1 (x)}\right]'=\Frac{\beta_1 (x)}{\beta_0 (x)}
\left[\Frac{\eta_0' (x)}{\sqrt{\eta_0 (x)^2+1}}-\Frac{\eta_1' (x)}{\sqrt{\eta_1 (x)^2+1}}\right],
$$
where $\eta_0(x)$ is given by (\ref{otrosc}) taking $\alpha=0$ and, using Eq. (\ref{eta1I}), 
$$
\eta_1(x)=f(x)\eta_0 (x),\,f(x)=1+\Frac{1}{2\sqrt{(\nu-1/2)^2 +x^2}}.
$$
Now, because $\eta_0 ' (x)<0$ and $\eta_1'(x)<0$ if $\nu>1/2$ then $\phi_1 '(x) >0$
is equivalent to 
\begin{equation}
\label{te1}
1+\eta_1(x)^2<\left(\Frac{\eta_1 '(x)}{\eta_0 '(x)}\right)^2 (1+\eta_0 (x)^2).
\end{equation}
Now, using $\eta_1 (x)=f(x)\eta_0 (x)$, the inequality (\ref{te1}) becomes
$$
1+\eta_0(x)^2 f(x)^2<\left(f(x)+\Frac{\eta_0 (x)}{\eta_0 '(x)}f'(x)\right)^2(1+\eta_0 (x)^2)
$$ 
and this inequality clearly holds because $\eta_0 (x)f'(x)/\eta_0 '(x)>0$ and $f(x)>1$.

\vspace*{0.2cm}
\noindent
{\bf 2. {\boldmath $\alpha=2$:}}

First we prove the case $\nu=0$, which reads
$$
\Frac{I_0(x)}{I_1(x)}>\Frac{x}{\delta+\sqrt{\delta^2+x^2}},\,\delta=-\Frac{1}{2}\left(1-\Frac{1}{\sqrt{1+4x^2}}\right).
$$
Now, because 
$\delta\in (-1/2,0)$ if $x>0$, we have $\Frac{x}{\delta+\sqrt{\delta^2+x^2}}<\Frac{x}{-1/2+\sqrt{\frac14+x^2}}$ and therefore
it is enough to prove that
$$
\Frac{I_0(x)}{I_1(x)}>\Frac{x}{-1/2+\sqrt{\frac14+x^2}}=\Frac{1}{x}\left(\frac12+\sqrt{\frac14+x^2}\right),
$$
which is true on account of Theorem \ref{boundeb} (setting $\nu=1$ in the first inequality).

Now we consider $\nu>0$. Using (\ref{serh1}) we see that $h_1 (0^+)>0$, $h_1^{\prime} (0^{+})<0$. Therefore, if we can prove that
$\phi_1 '(x)<0$, this will show that $h_1(x)>\phi_1(x)$ and then $I_{\nu}(x)/I_{\nu-1}(x)=x^2 h_0(x)>x^2 \beta_1 (x)=B_2 (\nu,x)$.

For $\alpha=2$ we have, for $i=0,1$,
$$
\begin{array}{l}
\eta_i(x)=\lambda x q_i(x),\, q_0(x)=1,\,q_{1}(x)=f(x),\\
\\f(x)=\gamma +\Frac{1}{2\sqrt{\lambda^2 +x^2}},
\lambda =\nu +\Frac{1}{2},\, \gamma=\Frac{\nu -\frac12}{\nu +\frac12}. 
\end{array}
$$
and
$$
\begin{array}{l}
\beta_i (x)=(\eta_i (x)+\sqrt{\eta_i (x)^2+x^4})^{-1}.
\end{array}
$$
Now, we compute the derivative $\phi_1^{\prime}(x)$:
$$
\phi_{1}'(x)=\Frac{\beta_1 (x)}{\beta_0 (x)}\left(\Frac{\beta_1^{\prime} (x)}{\beta_1 (x)}-\Frac{\beta_0^{\prime} (x)}{\beta_0 (x)}\right),
$$
with 
$$
-\Frac{\beta_i^{\prime}(x)}{\beta_i (x)}=\Frac{\eta_i^{\prime}(x)+2x^3 \beta_i (x)}{\sqrt{\eta_i (x)^2+x^4}}=\Frac{1}{x}
\left[2+\Frac{\eta_i^{\prime}(x)-2\lambda q_i (x)}{\sqrt{(\lambda q_i (x))^2+x^2}}\right].
$$
With this, $\phi_1^{\prime}(x)<0$ is equivalent to
$$
\Frac{f(x)-xf'(x)}{\sqrt{(\lambda f(x))^2+x^2}}<\Frac{1}{\sqrt{\lambda^2+x^2}}
$$
Of course, this inequality holds of $f(x)-xf'(x)<0$. Let us then assume that $f(x)-xf'(x)>0$ and let us prove that it holds
for $\nu\ge 0$ in this case. 

After assuming that $f(x)-xf'(x)>0$ we take squares and then $(f(x)-xf'(x))^2 (\lambda^2+x^2)<(\lambda f(x))^2+x^2$, and after some
simplification 
$$
\Frac{1}{x^2+\lambda^2}-(1-\gamma^2)-\Frac{\lambda^2}{4(\lambda^2+x^2)^2}<-2\gamma \Frac{1}{\sqrt{\lambda^2+x^2}}.
$$
We omit the last term in the left-hand side and consider the following inequality 
(which, of course, implies the previous inequality):
$$
\Frac{1}{\lambda^2 +x^2}-(1-\gamma^2)<-\Frac{2\gamma}{\sqrt{\lambda^2 +x^2}}.
$$
Now, letting $\zeta=\sqrt{\lambda^2 +x^2}$ the previous inequality is equivalent to 
$$
(1-\gamma^2)\zeta^2-2\gamma\zeta-1=(1-\gamma^2)(\zeta-\lambda)(\zeta+1/(1+\gamma))>0
$$
which is true if $|\gamma|<1$ because $\zeta>\lambda$. Therefore it holds for $\nu\ge 0$.

\hfill $\square$.

\begin{remark}
As said in the proof, $B_{\alpha}(\nu,x)$ is decreasing as a function of $\alpha$. Therefore $B_{\alpha}(\lambda,x)$ are bounds
for $\alpha\ge 2$ and in particular $B_{+\infty}(\nu,x)$ is a lower bound, as we already know from Theorem \ref{boundeb} 
($B_{+\infty}(\nu,x)=
b_1(\nu,x)$). Similarly, $B_{1-2\nu}(\nu,x)=b_0(\nu,x)$ is a bound for $\nu\ge 1/2$, as was also shown in Theorem \ref{boundeb}. 
Then we see that the bounds in Theorem \ref{cotasric} are sharper than the bounds in Theorem \ref{boundeb}. If $\nu>1/2$ 
the upper bound in Theorem \ref{cotasric} is of course also sharper than the bound of Theorem \ref{masdeb}, which is 
$b_{-1}(\nu,x)=B_{-\infty}(\nu,x)$ (but this bound is valid for $\nu\ge 0$). 
\end{remark}

After iterating once the Riccati equation, it is natural to ask what is the result of the next iteration. Numerical experiments show
that, although the new approximations are sharper, particularly for large parameters, they are not bounds for any real $\alpha$.
However, as before, we can use the iteration of the recurrence relation (\ref{CF}) in order to try to improve the bounds. We obtain the following

\begin{theorem}
\label{iteric}
Let 
\begin{equation}
\begin{array}{l}
\widetilde{B}_{\alpha} (\nu,x)=\Frac{x}{\delta^{-}_{\alpha} (\nu,x)+
\sqrt{\delta^{+}_{\alpha} (\nu ,x)^2+x^2}},\\
\delta_{\alpha}^{\pm} (\nu,x)=(\nu\pm 1/2)\pm\Frac{\sigma}{2\sqrt{\sigma^2+x^2}},\,\sigma=\nu+(\alpha+1)/2,
\end{array}
\end{equation}
then
$$
\widetilde{B}_0 (\nu,x)<\Frac{I_{\nu}(x)}{I_{\nu-1}(x)}<\widetilde{B}_2 (\nu,x),\,\nu \ge 0
$$
\end{theorem}

The lower bound in this theorem turns out to be sharper than that from Theorem \ref{cotasric}; however, the same is
not true for the upper bound, and for large enough $x$ the bound in Theorem \ref{cotasric} is sharper. As we discuss next,
the iteration of the recurrence, differently from the iteration of the Riccati equation, does not improve the order
of approximation for large $x$.

\subsection{Comparison between bounds}

From the discussion after Theorems \ref{boundeb2} and \ref{iteric} we observe that the iteration of the recurrence
relation may in some cases lead to improved bounds but that for large $x$ some of these bounds are not improved.
This is not surprising because the iteration of the continued fraction (\ref{CF}) does not improve the sharpness of
the bounds as $x\rightarrow +\infty$, differently to the cases $x\rightarrow 0$ and $\nu\rightarrow +\infty$.

In order to see this, let us consider the sequences of bounds generated by (\ref{CF}) starting from the
lower and upper bounds $l_{\nu}^{(0)}(x)$, $u_{\nu}^{(0)}(x)$, where this initial bounds can be given by Theorem
\ref{boundeb} or Theorem \ref{cotasric}. The first iteration of the CF is given by Theorems \ref{boundeb2} and 
\ref{iteric}. We have the sequences of lower and upper bounds:
\begin{equation}
\label{seqb}
l_{\nu}^{(i+1)}(x)=\Frac{1}{\Frac{2\nu}{x}+u_{\nu+1}^{(i)}(x)},\,
u_{\nu}^{(i+1)}(x)=\Frac{1}{\Frac{2\nu}{x}+l_{\nu+1}^{(i)}(x)}.
\end{equation}
In order to see whether the bounds are sharper in each iteration we can check how the lower and upper bounds come closer
in each iteration. For measuring this, we consider $c_{\nu}^{(i)}(x)=u_{\nu}^{(i)}(x)/l_{\nu}^{(i)}(x)-1$. From
(\ref{seqb}) we get:
$$
c_{\nu}^{(i+1)}(x)=c_{\nu +1}^{(i)}(x)\Frac{1}{1+\Frac{2\nu}{x l_{\nu+1}^{(i)}(x)}}
$$
Now, we have that, either starting from Theorem \ref{boundeb} or Theorem \ref{cotasric}: 
$$l^{(i)}_{\nu}(x)=\Frac{x}{2\nu}(1+{\cal O}(\nu^{-1})),\nu\rightarrow +\infty,$$
and similarly as $x\rightarrow 0$,
while $l^{(i)}_{\nu}(x)=1+{\cal O}(x^{-1})$ as $x\rightarrow +\infty$. Therefore
\begin{equation}
c_{\nu}^{(i+1)}(x)\sim \Frac{x^2}{4\nu^2}c_{\nu +1}^{(i)}(x)
\end{equation}
as $x\rightarrow 0$ or $\nu\rightarrow +\infty$ while
\begin{equation}
c_{\nu}^{(i+1)}(x)= c_{\nu +1}^{(i)}(x) \left(1+{\cal O}(x^{-1})\right),\,x\rightarrow +\infty
\end{equation}

This also explains why the new bounds of Theorem \ref{cotasric} have an intrinsic advantage over the previous bounds
(Theorem \ref{boundeb}) or their first iteration 
(Theorem \ref{boundeb2}): no matter how many times the CF is iterated, for a fixed number of iterations 
they are sharper for large enough $x$.

Considering the bounds from Theorem \ref{boundeb} we have
$$
\begin{array}{l}
c_{\nu}^{(0)}(x)=\Frac{1}{2\nu-1}(1+{\cal O}(x)),\,\nu>1/2,\,x\rightarrow 0\\ 
\\
c_{\nu}^{(0)}(x)=\Frac{1}{2\nu}(1+{\cal O}(\nu^{-1})),\,\nu\rightarrow +\infty\\
\\
c_{\nu}^{(0)}(x)=\Frac{1}{2 x}(1+{\cal O}(x^{-1})),\,x\rightarrow +\infty .
\end{array}
$$
On the the other hand denoting by $\tilde{c}_{\nu}^{(0)}(x)$ the same quantities, but using the bounds in Theorem \ref{cotasric}
we have
$$
\begin{array}{l}
{\tilde c}_{\nu}^{(0)}(x)=\Frac{8x^2}{(4\nu^2-1)^2}(1+{\cal O}(x)),\,\nu>1/2,\,x\rightarrow 0\\ 
\\
{\tilde c}_{\nu}^{(0)}(x)=\Frac{x^2}{2\nu^4}(1+{\cal O}(\nu^{-1}),\,\nu\rightarrow +\infty\\
\\
{\tilde c}_{\nu}^{(0)}(x)=\Frac{1}{2 x^2}(1+{\cal O}(x^{-1})),\,x\rightarrow +\infty .
\end{array}
$$

As we know, the bounds in Theorem \ref{cotasric} are sharper than those of Thm. \ref{boundeb}. The previous estimations
confirm this fact in the three different limits and, furthermore, show
that the bounds from \ref{cotasric} are, for large enough $x$, sharper than the iterated bounds starting from Thm. \ref{boundeb}.
Indeed, the CF improves the sharpness of the bounds as $x\rightarrow 0$ and $\nu\rightarrow +\infty$, and we have
$$
c_{\nu}^{(i)}(x)=\Frac{1}{2\nu}\left(\Frac{x}{2\nu}\right)^{2i}\left(1+{\cal O}(x,\nu^{-1})\right),\,
{\tilde c}_{\nu}^{(i)}(x)=\Frac{1}{2\nu^2}\left(\Frac{x}{2\nu}\right)^{2i+2}\left(1+{\cal O}(x,\nu^{-1})\right)
$$
but as $x\rightarrow +\infty$
\begin{equation}
\label{xlarge}
c_{\nu}^{(i)}(x)=\Frac{1}{2x}\left(1+{\cal O}(x^{-1})\right),\,
{\tilde c}_{\nu}^{(i)}(x)=\Frac{1}{2x^2}\left(1+{\cal O}(x^{-1})\right).
\end{equation}

The bounds of Theorem \ref{cotasric} can be used for improving the computation of the ratio $I_{\nu}(x)/I_{\nu-1}(x)$
by means of the continued fraction arising from (\ref{CF}): using Theorem \ref{cotasric} as tail estimation 
improves the convergence for large $x$. Notice that the approximants of the continued fraction expansion for the ratio come from
iterating (\ref{CF}) using
$$
H_{\nu}^{(i+1)}(x)=\Frac{1}{\Frac{2\nu}{x}+H_{\nu+1}^{(i)}(x)},\, H_{\nu}^{(0)}(x)\equiv 0,
$$
where $H_{\nu}^{(2i-1)}(x)$, $i=1,2,\ldots$ are upper bounds and $H_{\nu}^{(2i)}(x)$ are lower bounds for $I_{\nu}(x)/I_{\nu-1}(x)$.
We have in this case
$$
\Frac{H_{\nu}^{(2i-1)}(x)}{H_{\nu}^{(2i)}(x)}-1=\Frac{1}{(2i)^2 (\nu+i)(\nu+i-1)}x^2\left(1+{\cal O}(x^{-2})\right),
$$
which has the wrong asympotics as $x$ becomes 
large. A first improvement comes by using Thm. \ref{boundeb} which gives $c_{\nu}^{(i)}={\cal O}(x^{-1})$ and
 a further improvement comes from Thm. \ref{cotasric} which gives $c_{\nu}^{(i)}={\cal O}(x^{-2})$ (see eq. (\ref{xlarge})).

\section{Bounds for the ratio $K_{\nu}(x)/K_{\nu-1}(x)$}

We consider a similar analysis for the ratio $K_{\nu}(x)/K_{\nu -1}(x)$. We start with 
\begin{equation}
\label{dk1}
h (x)=K_{\nu-1}(x)/K_{\nu}(x), 
\end{equation}
and 
we take 
\begin{equation}
\label{dk2}
h_0(x)=x^{-\alpha}h (x), 
\end{equation}
which satisfies
\begin{equation}
\label{ricaa2}
h_0 '(x)=-x^{-\alpha}+\Frac{2\tau}{x}h_0 (x)+x^{\alpha} h_0 (x)^2,\,\tau=\nu-\frac12 (\alpha +1).
\end{equation}
We have 
\begin{equation}
\label{dk3}
\gamma (x)=x^{2\alpha},\,\eta_0 (x) =\tau x^{\alpha-1}
\end{equation}
and the potential bound for 
$h_0(x)$ is
\begin{equation}
\label{bb2}
\beta_0 (x)=\phi_0(x)=\Frac{1}{\eta_0(x)+\sqrt{\eta_0 (x)^2+\phi (x)}}=\Frac{x^{1-\alpha}}{\tau+\sqrt{\tau^2 +x^2}}.
\end{equation}
The uniparametric family  (with parameter $\alpha$) of possible bounds for $h(x)$ is then
\begin{equation}
\label{ba2}
d_{\alpha} (\nu,x)=\Frac{x}{\tau+\sqrt{\tau^2 +x^2}},\,\tau=\nu-\frac12 (\alpha +1).
\end{equation}
As it is known, $d_{0} (\nu,x)$ is an upper bound for $h(x)$ when $\nu\ge 1/2$, as was proved in
\cite{Segura:2011:BFR} using the Riccati equation technique; 
$d_{1} (\nu, x)$ is also known to be a lower bound, as can be proved, for instance, using a Tur\'an-type inequality
(see \cite{Laforgia:2010:SIF,Segura:2011:BFR}). Here we prove both results using Theorem \ref{esode} and establish that
these are the best possible bounds of this form.

\begin{theorem}
\label{K0}
$$
\Frac{K_{\nu-1}(x)}{K_{\nu}(x)} \ge d_{0}(\nu,x)=\Frac{x}{\nu- 1/2+\sqrt{(\nu-1/2)^2+x^2}},\,\nu\ge 1/2,
$$
where the equality only holds for $\nu=1/2$.
$$
\Frac{K_{\nu-1}(x)}{K_{\nu}(x)} < d_{1}(\nu,x )=\Frac{x}{\nu- 1+\sqrt{(\nu -1)^2+x^2}},\,\nu\in {\mathbb R}.
$$ 
These bounds are the best possible of the form $d_{\alpha} (\nu,x)$, $\alpha\in {\mathbb R}$, in their range of
validity.
\end{theorem}

\noindent {\it Proof.} The bounds correspond to the values $\alpha=0,1$, 

First we check that $d_{\alpha}(\nu,x)$ are not bounds for $\alpha\in (0,1)$. To prove this, it is enough to consider
$\nu>1$. Using that $K_{\mu}(x)\sim \frac12 \Gamma (\mu) (x/2)^{-\mu}$ for $\nu>0$ as $x\rightarrow 0^+$ 
and the asymptotic expansion
\cite[10.40.2]{Olver:2010:BF}, we have:
\begin{equation}
\label{hh2}
\begin{array}{l}
h_0(x)=\Frac{x^{1-\alpha}}{2(\nu-1)}\left(1+ o(1)\right),\,\nu>1,\\
h_0(x)=x^{-\alpha}\left(1-(\nu-1/2)x^{-1}+{\cal O}(x^{-2})\right).
\end{array}
\end{equation}
Therefore $h_0(0^+)h_0(+\infty)<0$ if $\alpha\in (0,1)$; this implies that the graph of $h_0(x)$ must cross
the graph of the characteristic root $\phi_0(x)$, which therefore can not be a bound for $h_0 (x)$.

Next, we prove the bounds for $\alpha=0$ (upper bound) and $\alpha=1$ (lower bound)

\begin{enumerate}
\item{$\alpha=0$:} the case $\nu=1/2$ is obvious because $K_{1/2}(x)=K_{-1/2}(x)$. For $\nu>1/2$ 
we have that $\phi_0 (x)=x/(\tau+\sqrt{\tau^2+x^2})$, $\tau=\nu -1/2$, is strictly increasing\footnote{Observe that 
for this to hold we need $\tau>0$, that is, $\nu>1/2$}; in
addition, (\ref{hh2}) shows that $h_0(+\infty)>0$ and $h_0 '(+\infty)>0$ 
(and then $\phi_0 ' (+\infty)h_0 '(+\infty)>0$). Applying Theorem \ref{esode} for
the case $C(x)>0$, we have $h_0(x)=K_{\nu -1}(x)/K_{\nu}(x)>\phi_0 (x)$.
\item{$\alpha=1$:} we have that $\phi_0 (x)=1/(\nu-1+\sqrt{(\nu-1)^2 +x^2})$
is strictly decreasing as a function of $x$; in
addition, (\ref{h0i}) shows that $h_0 '(+\infty)<0$, and then $\phi_0 ' (+\infty)h_0 '(+\infty)>0$. Applying \ref{esode}, 
we have $h_0(x)= x^{-1} K_{\nu -1}(x)/K_{\nu}(x)<\phi_0 (x)$.

\end{enumerate}

Finally, we notice that $d_{\alpha} (\nu, x)$ is increasing as
a function of $\alpha$, and then, because $d_1 (\nu,x)$ is an upper bound, $d_{\alpha} (\nu,x)$ are upper bounds for 
any $\alpha\ge 1$; similarly, because $d_0 (\nu,x)$ is a lower bound, $d_{\alpha} (\nu,x)$ are lower bounds for any 
$\alpha\le 0$. The sharpest bounds correspond to $\alpha=0,1$. 

\hfill $\square$.

Similarly as for the $I$ Bessel function, Theorem \ref{esode} can not be used to prove that any of the values $\beta_0 (x)$ is a bound for values
$\alpha \in (-1,0)\cup(0,1)$, because $\phi_0 (x)=\beta_0 (x)$ is not monotonic in this case. 
However, from the monotonicity of $d_\alpha (\nu,x)$ we concluded
that for $\alpha\in (-1,0)$ we have lower bounds for $K_{\nu-1}(x)/K_{\nu}(x)$ 
at least for $\nu\ge 1/2$; in fact, these bounds, although less sharp than 
$d_0(\nu,x)$ for positive $\nu$, extend the range of validity with respect. For $\alpha=-1$ we obtain the following result:

\begin{theorem}
\label{otrok}
$$
\Frac{K_{\nu -1}(x)}{K_{\nu}(x)}>d_{-1} (\nu,x)=\Frac{x}{\nu+\sqrt{\nu^2+x^2}},\,\nu\in {\mathbb R},
$$
\end{theorem}

\noindent {\it Proof.} For any real $\nu$ we have that $\phi_0(x)$ is increasing for $\alpha=-1$ 
and $h_0 '(0^+)>0$, and
Theorem \ref{esode} implies that $h_0(x)=xK_{\nu-1}(x)/K_{\nu}(x)>\phi_0(x)=x^2/(\nu+\sqrt{\nu^2+x^2})$.

\hfill $\square$.

\begin{remark}
Theorem \ref{otrok} is in fact equivalent to the upper bound in Theorem \ref{K0}. Indeed, using that 
$K_{\mu}(x)=K_{-\mu}(x)$ we have that
$$
\Frac{K_{\nu-1}(x)}{K_{\nu}(x)}=
\left(\Frac{K_{-\nu}(x)}{K_{1-\nu}(x)}\right)^{-1}>
\Frac{1}{d_1 (1-\nu,x)}=\Frac{x}{\nu+\sqrt{\nu^2+x^2}}
$$
\end{remark}

\vspace*{0.3cm}
As for the case of the Bessel function of the first kind, it is also possible to obtain new bounds
by iterating the recurrence relation, as was done in \cite{Segura:2011:BFR}. In this case we have
\begin{equation}
\label{cfk}
\Frac{K_{\nu -1}(x)}{K_{\nu}(x)}=\Frac{1}{\Frac{2(\nu -1)}{x}+\Frac{K_{\nu -2}(x)}{K_{\nu -1}(x)}}.
\end{equation}
Differently from the previous case, these bounds from the iteration of (\ref{cfk}) are of more restricted validity 
with respect to $\nu$ as more
iterations are considered. We will not describe this type bounds, neither starting from the bounds
from Theorems \ref{K0} and \ref{otrok}
nor with the improved bounds we discuss next. 

\subsection{Bounds from the iteration of the Riccati equation}

We start from (\ref{dk1}), (\ref{dk2}) and (\ref{dk3}). The first iteration (following the notation of 
section \ref{iterating}) gives
\begin{equation}
\label{eta1k}
\eta_1 (x)=x^{\alpha-1}\left(\nu-\frac12-\Frac{\tau}{2\sqrt{\tau^2+x^2}}\right),\,\tau=\nu-\Frac{\alpha+1}{2},
\end{equation}
and then the potential bounds for $K_{\nu -1}(x)/K_{\nu}(x)$ are of the form
\begin{equation}
\label{bok}
D_{\alpha}(\nu,x)=\Frac{x}{\varphi_{\alpha} (\nu,x)+
\sqrt{\varphi_{\alpha} (\nu,x)^2+x^2}},
\end{equation}
with
$$\varphi_{\alpha} (\nu,x)=(\nu-1/2)-\Frac{\tau}{2\sqrt{\tau^2+x^2}}.$$ 
We will apply again Theorem \ref{esode} for checking whether the function $h_1 (x)=h_0(x)/\beta_0(x)$ 
is bounded by $\phi_1(x)=\beta_1 (x)/\beta_0(x)$. For proving this, we will need the expansion
\begin{equation}
\label{h1exk}
h_1(x)=1-\Frac{\alpha}{2}x^{-1}+\Frac{1}{2}\left(\nu -\frac12 +\frac14 \alpha^2\right) x^{-2}+{\cal O}(x^{-3})
\end{equation} 

In particular, we prove that for $\alpha=0$ we obtain an new upper bound for $h(x)$, which improves the 
upper bound of Theorem \ref{K0}, but that it is not possible to improve the lower bound in its range of validity.

\begin{theorem}
The following bounds hold for $\nu\ge 1/2$
\begin{equation}
\begin{array}{l}
D_{2\nu-1}(\nu,x)\le \Frac{K_{\nu-1}(x)}{K_{\nu}(x)}\le D_0(\nu,x)
\end{array}
\end{equation}
where the equalities only hold for $\nu=1/2$. No bound $D_{\alpha}(\nu,x)$
for fixed $\alpha\neq 0$ exists ($\alpha$ not depending on $\nu$) sharper than 
than $D_{2\nu-1}(\nu,x)$ or $D_0 (\nu,x)$ for $\nu\ge 1/2$.
\end{theorem}

\noindent {\it Proof.} 
The lower bound was already proved in Theorem \ref{K0} ($D_{2\nu-1}(\nu,x)=d_0 (\nu,x)$). 
Once we prove that $D_0(\nu,x)$ is an upper bound, with the equality holding for $\nu=1/2$,
it will be obvious that there is no other bound of the tyoe $D_{\alpha}(\nu,x)$ sharper than 
$D_{2\nu-1}(\nu,x)$ or $D_0 (\nu,x)$ for $\nu\ge 1/2$; clearly, such bound can not be
as sharp as $D_{2\nu-1}(\nu,x)$ and $D_0 (\nu,x)$ as $\nu\rightarrow 1/2$.

Now we prove that the bound obtained
from the next iteration of the Riccati equation for $\alpha=0$ ($D_0(\nu,x)$) is an upper bound. 
For proving this, that is, for proving that 
$h_0(x)=K_{\nu-1}(x)/K_{\nu}(x)<\beta_1 (x)$
for $\nu >1/2$ (the case $\nu=1/2$ is trivial),
 we only need to compare the monotonicity of $\phi_1(x)=\beta_1 (x)/\beta_0 (x)$ against
the monotonicity of $h_1(x)=h_0(x)/\beta_0(x)$ for large $x$.

With respect to the behavior of $h_1 (x)$ we have,
using the asymptotic expansion (\ref{h1exk} that
$h_1(x)>0$, $h_1'(x)<0$ for $\nu>1/2$ and $x$ large enough.

Next we prove that $\phi_1 (x)=\beta_1 (x)/\beta_0 (x)$ is decreasing for all $x>0$ if $\nu>1/2$. We have
$$
\phi_1'(x)=\Frac{\beta_1 (x)}{\beta_0 (x)}
\left[\Frac{\eta_0' (x)}{\sqrt{\eta_0 (x)^2+1}}-\Frac{\eta_1' (x)}{\sqrt{\eta_1 (x)^2+1}}\right]
$$
where (taking $\alpha=0$ in (\ref{dk3}) and (\ref{eta1k})) 
$\eta_1 (x)=\eta_0 (x)f(x)$ with $f(x)=1-\Frac{1}{2\sqrt{(\nu-1/2)^2+x^2}}<1$. Two cases are possible: $\eta_1 ' (x)\ge 0$ or 
$\eta_1 '(x)<0$. In the fist case, it is obvious that $\phi_1' (x)<0$ because $\eta_0'(x)$ is always negative. In the second case,
we have that $\phi_1 '(x)<0$ is equivalent to
\begin{equation}
\label{lac}
1+\eta_0 (x)^2 f(x)^2 > \left(f(x)+\Frac{\eta_0 (x)}{\eta_0 ' (x)}f'(x)\right)^2 (1+\eta_0 (x)^2)
\end{equation}
Now, we observe that the term under parenthesis, which is equal to $\eta_1 '(x)/\eta_0 '(x)$, satisfies
\begin{equation}
\label{desm}
0<\Frac{\eta_1 '(x)}{\eta_0 '(x)}=f(x)+\Frac{\eta_0 (x)}{\eta_0 ' (x)}f'(x)< f(x)<1.
\end{equation}
The left-hand side equality is obvious because $\eta_0 '(x)<0$ and we are assuming that $\eta_1 '(x)<0$; that 
$\eta_1 '(x)/\eta_0 '(x)<f(x)$ is also obvious. And using (\ref{desm}), (\ref{lac}) is proved.

Therefore, for the case $\nu>1/2$ we have proved that $h_1 (x)>0$, $h_1 '(x)<0$ for large $x$ and that the characteristic root 
$\phi_1 (x)$ is decreasing for $x>0$. Using similar arguments as those considered
 in the proof of Lemma 2, this implies that $h_1 (x)>0$, $h_1'(x)<0$ 
for all $x$ and that $h_1 (x)<\phi_1(x)$. Therefore $h_0(x)<\phi_1(x)\phi_0(x)=\beta_1(x)$ and the theorem is proved. \hfill $\square$

\begin{remark}
$D_{\alpha}(\nu,x)$ is decreasing as a function of $\alpha$. Then, because $D_{2\nu-1}(\nu,x)$ is a lower bound,
$D_{\alpha}(\nu,x)$ are also lower bounds (but weaker) when $\alpha\ge 2\nu-1$. Similarly $D_{\alpha}(\nu,x)$ are 
upper bounds for $\alpha\le 0$ (again, weaker). We already proved that $D_{\mp \infty}(\nu,x)$ are bounds 
(Theorems \ref{K0} and \ref{otrok}),
which are necessarily less sharp than  $D_0(\nu,x)$ and $D_{2\nu-1}(\nu,x)$ respectively.
\end{remark}

\section{Concluding remark}

Similar ideas can potentially
establish new sharp bounds for other functions. For instance, some of the bounds in \cite{Segura:2014:MPA}
can be reobtained using these ideas and the optimality of some of these bounds can be proven easily.

Both Bessel and incomplete gamma functions can be expressed in terms of confluent hypergeometric functions,
and it is tempting to consider this more general case and to study when these technieques can be used
for building improved bounds. More generally, hypergeometric function ratios, and particularly confluent ratios, satisfy Riccati
equations from which sharp bounds can be extracted (a recent example of this can be found in 
the appendix of reference \cite{Segura:2015:SBF}). We expect that the ideas used in the present paper can also
be used to obtain new and improved bounds in this more general case.

\section{Acknowledgements}
JS acknowledges financial support from Ministerio de Econom\'{\i}a y Competitividad, project MTM2012-34787. DRA
acknowledges a predoctoral research contract (FPI) from Ministerio de Econom\'{\i}a y Competitividad 
(reference BES-2013-064743).

\bibliographystyle{plain}

\bibliography{bessbarx}

\end{document}